\documentclass[a4paper,11pt,fullpage]{article}
\usepackage[english]{babel}

\usepackage[left=1in,right=1in,top=1in,bottom=1in]{geometry}

\usepackage{amsmath}
\usepackage{amsthm}
\usepackage{amsfonts}
\usepackage{amssymb}
\usepackage{amsxtra}
\usepackage{latexsym}
\usepackage{ifthen}
\usepackage{cite}
\usepackage{esint}
\usepackage[cp1250]{inputenc}

\usepackage{epic}
\usepackage{eepic}
\usepackage{stmaryrd}

\usepackage{cancel}
\usepackage[normalem]{ulem}
\usepackage[dvipsnames]{xcolor}%

\DeclareMathOperator*{\esssup}{ess\,sup}

\numberwithin{equation}{section}
\newtheorem{theorem}{Theorem}[section]
\newtheorem{lemma}{Lemma}[section]
\newtheorem{remark}{Remark}[section]

\def\XXint#1#2#3{{\setbox0=\hbox{$#1{#2#3}{\int}$}
     \vcenter{\hbox{$#2#3$}}\kern-.5\wd0}}

\begin{document}

\title{On the local boundedness of weak solutions to doubly nonlinear equations with variable exponents
}

\author{Simone Ciani, Mariia O. Savchenko, Igor I. Skrypnik, Yevhen S.  Zozulia 
}

\date{}

  \maketitle

\begin{abstract} 
\noindent 
We establish local boundedness, together with quantitative estimates, for non-negative weak sub-solutions to a class of doubly nonlinear parabolic equations whose prototype is
\[
(u^q)_t-\operatorname{div}
\Big(|Du|^{p(x,t)-2}Du\Big)=0,
\qquad p(x,t)>1,\quad q>0.
\]
The results are obtained under a logarithmic continuity condition on the variable exponent $p(x,t)$ and suitable relations between $q$ and the local bounds of $p(x,t)$. In certain ranges of the exponents, an additional local integrability assumption is required.

\vskip0.2cm 


\end{abstract}

\pagestyle{myheadings} \thispagestyle{plain}
\markboth{  }
{Doubly nonlinear parabolic equations... }

\section{Introduction  and Main Results}\label{Introduction}

Nonlinear parabolic equations with nonstandard growth conditions arise naturally in various models of continuum mechanics and have attracted considerable attention in recent decades.
An important class of such equations is associated with the mathematical modeling of electrorheological fluids, whose mechanical properties depend on the strength of an applied electric field; see \cite{Pao,Ruzicka-book,exp1}. The corresponding constitutive laws lead to differential operators with variable exponent growth, reflecting the nonhomogeneous rheological behaviour of these materials. A typical  model is the $p(x)$-Laplace equation
\[
-\operatorname{div}\bigl(|Du|^{p(x)-2}Du\bigr)=0,
\]
whose mathematical analysis is closely related to the theory of Lebesgue and Sobolev spaces with variable exponents; see \cite{DHHR-book}.  Such operators generalize the classical $p$-Laplacian and give rise to a rich mathematical structure due to the loss of homogeneity and uniform ellipticity.
In evolutionary problems, the exponent may depend on both spatial and temporal variables, leading to parabolic equations of the form
\[
u_t-\operatorname{div}
\bigl(|Du|^{p(x,t)-2}Du\bigr)=0.
\]
The existence and regularity properties of weak solutions to such equations have been extensively investigated. Antontsev and Shmarev \cite{AntShm2007,AntShm2009} established existence results and boundedness estimates for weak solutions to more general anisotropic parabolic equations with variable exponent growth. Xu and Chen \cite{XuChe} proved the local boundedness of weak solutions and, under the additional assumption that the lower bound of $p(x,t)$ exceeds $2$, established their interior H\"older continuity. Subsequently, Alkhutov and Zhikov \cite{AlkZik} obtained H\"older continuity results for solutions to parabolic equations with variable exponents of nonlinearity. Wang \cite{Wan} established intrinsic Harnack inequalities for parabolic equations with variable exponents. B\"ogelein and Duzaar \cite{BogDuz2012} proved the local H\"older continuity of the spatial gradient of weak solutions to parabolic $p(x,t)$-Laplacian systems under suitable H\"older continuity assumptions on the exponent. Yao \cite{Yao15} extended the gradient regularity theory by establishing local H\"older continuity of the spatial gradient of weak solutions to general nonhomogeneous parabolic equations with $p(x,t)$-growth. Furthermore, Baroni and B\"ogelein \cite{BarBog14} obtained local Calder\'on--Zygmund estimates for parabolic $p(x,t)$-Laplacian systems under a logarithmic continuity assumption on $p(x,t)$, showing that the spatial gradient inherits higher integrability properties from the prescribed data.

In this work, we investigate the local boundedness of weak solutions to a class of doubly nonlinear parabolic equations with variable exponent growth, whose prototype is
\[
(u^q)_t-\operatorname{div}\Big(|Du|^{p(x,t)-2}Du\Big)=0,
\qquad p(x,t)>1,\quad q>0.
\]
For the constant exponent case $p(x,t)\equiv p>1$, equations of this type are well known, have numerous applications, and have been studied extensively for several decades (see, for instance, \cite{Kalashnikov,Raviart,Kalashnikov-survey} and the references therein). We assume that the exponent $p(x,t)$ is bounded and satisfies the logarithmic continuity condition (see \eqref{eq1.3}), which is weaker than H\"older continuity. For related models of electrorheological fluids with nonstandard continuity assumptions on the exponent, we refer to \cite{KR}. Under suitable relations between $q$ and the local bounds of $p(x,t)$, we establish quantitative local boundedness estimates, extending the regularity theory for parabolic equations with variable exponent growth to the doubly nonlinear setting.

\subsection{Notations}

\noindent We provide some notations that will be used in the sequel, more precisely,
for fixed $(x_0, t_0)\in \Omega_T$, and $r, \eta >0$ we define the following  cylinders 
$$Q_{r, \eta}(x_0, t_0):= B_r(x_0)\times I_\eta(t_0), \quad I_\eta(t_0):=(t_0-\eta, t_0).$$
For a positive continuous function $p(x, t)$ we denote also 
$$ p^+_{r, o}:=\sup\limits_{Q_{r, r}(x_0, t_0)}p(x, t),\qquad p^-_{r, o}:=\inf\limits_{Q_{r, r}(x_0, t_0)} p(x, t).$$
We define the Sobolev spaces of  integrable functions with variable exponents as
$$W (\Omega):=\big\{u\in W^{1,1} (\Omega): \int\limits_{\Omega}|D u|^{p(x, t)} dx < \infty, \quad \forall\,\, t\in(0, T)\big\},$$
$$L (0, T; W (\Omega)):=\big\{u\in L^1 (0, T ,{W^{1,1} (\Omega)}): \int\limits_{0}^{T}\int\limits_{\Omega}|D u|^{p(x, t)}\,dx dt < \infty\big\},$$
and the respective spaces of functions with zero boundary data
$$W_0(\Omega):= W^{1,1}_0(\Omega) \cap W (\Omega),$$
$$L (0, T; W_0(\Omega)):=L^1 (0, T; W^{1,1}_0(\Omega)) \cap L (0, T; W (\Omega)).$$
For any $a$, $b \in \mathbb{R}_+$ we also use the notation
$a \land b:=\min(a, b)$, $a \vee b:=\max (a, b)$.

\subsection{Main Results}
\noindent Let $\Omega$ be a bounded domain in $\mathbb{R}^{N}$, $T>0$, $\Omega_{T}:= \Omega \times (0, T)$,
we study  non-negative (sub) super-solutions to the equation
\begin{equation}\label{eq1.1}
(u^q)_t- div\mathbf{A}(x, t, u, D u)=0, \quad (x, t)\in \Omega_{T}.
\end{equation}
Throughout the paper we suppose that the functions $\mathbf{A}:\Omega_{T}\times \mathbb{R}_+\times \mathbb{R}^{N} \rightarrow \mathbb{R}^{N}$ are Carath\'eodory functions\footnote{Which is,   
$\mathbf{A}(\cdot, \cdot, u, \xi)$ are Lebesgue measurable for all  $u\in \mathbb{R}_+$, $\xi \in \mathbb{R}^{N}$,
and $\mathbf{A}(x, t, \cdot,  \cdot)$ are continuous for almost all $(x, t)\in \Omega_{T}$.}.
We also assume that the following structure conditions are satisfied
\begin{equation}\label{eq1.2}
\begin{cases}
\mathbf{A}(x, t, u, D u) D u \geqslant  K_{1}\,| D u|^{p(x, t)},\qquad \qquad 1<p_1\leqslant p(x, t)\leqslant p_2<\infty,\\
|\mathbf{A}(x, t, u,  D u)| \leqslant  K_{2}\,| D u|^{p(x, t)-1},
\end{cases}
\end{equation}
where $K_{1}$, $K_{2}$ are positive constants and the function $p(x, t): \Omega_{T} \rightarrow \mathbb{R}_{+}$
satisfies the so-called {\it logarithmic} condition: there exists a constant $K_3>0$ such 
\begin{equation}\label{eq1.3}
|p(x, t)-p(y, \tau)|\leqslant \frac{K_3}{\log\frac{1}{|x-y|+|t-\tau|}}\, \, ,
\end{equation}
for any $(x, t)$, $(y, \tau)\in \Omega_T$ such that $|x-y|+|t-\tau|\leqslant 1/2$. In what follows, the numbers $\{K_1, K_2, K_3, q, p_1, p_2, N\}$ are further referred to as {\it the data} and  we will write $\gamma$ as a generic positive constant that can be quantitatively determined a priori only in terms of {\it the data} if it depends on $\{K_1, K_2, K_3, q, p_1, p_2, N\}$, and that can change from line to line.
Postponing for the moment the formal definitions of weak solution, we will proceed to present the main results.

\begin{theorem}\label{th1.1}
Let $u$ be a non-negative, local, weak sub-solution to \eqref{eq1.1}-\eqref{eq1.3} in $\Omega_T$. Fix $(x_0, t_0)\in \Omega_T$
such that $Q_{\rho, \theta}(x_0, t_0)\subset Q_{\rho, \rho}(x_0, t_0)\subset Q_{8\rho, 8 \rho}(x_0, t_0)\subset \Omega_T$ and let
\begin{equation}\label{eq1.4}
q\leqslant p^-_{\rho, o}-1+(q+1)\frac{p^-_{\rho, o}}{N},
\end{equation}
then $u$ is bounded in $Q_{\rho, \theta}(x_0, t_0)$, provided that $\rho$ is small enough. Particularly, if
\begin{equation}\label{eq1.5}
q < p^-_{\rho, o}-1+(q+1)\frac{p^-_{\rho, o}}{N},
\end{equation}
then there exists $C>0$ depending only on the data such that if $p^+_{\rho, o}<q+1$, then 
\begin{multline}\label{eq1.6}
C^{-1} \sup\limits_{Q_{\rho, \theta}(x_0, t_0)} u\leqslant  \bigg(\frac{1}{\theta^{1+\frac{N}{p^-_{\rho, o}}}}
\iint\limits_{Q_{2\rho, 2\theta}(x_0, t_0)} u^{q+1}\,dx dt\bigg)^{\frac{p^-_{\rho, o}}{p^-_{\rho, o}(q+1+N)-(q+1)N}}+ \Big(\dfrac{\theta}{\rho^{p^+_{\rho, o}}}\Big)^{\frac{1}{q+1-p^+_{\rho, o}}} \, +\rho,
\end{multline} while if $q+1 < p^+_{\rho, o}$, then
\begin{multline}\label{eq1.7}
C^{-1}\sup\limits_{Q_{\rho, \theta}(x_0, t_0)} u\leqslant  \bigg(\frac{1}{\rho^{p^+_{\rho, o}(1+\tfrac{N}{p^-_{\rho, o}})}}
\iint\limits_{Q_{\rho, \theta}(x_0, t_0)} u^{p^+_{\rho, o}}\,dx dt\bigg)^{\frac{p^-_{\rho, o}}{p^-_{\rho, o}(q+1+N)-p^+_{\rho, o}N}}+ \Big(\dfrac{\rho^{p^+_{\rho, o}}}{\theta}\Big)^{\frac{1}{p^+_{\rho, o}-q-1}} + \rho,
\end{multline}
and when $p^+_{\rho, o}=q+1$ the estimate
\begin{multline}\label{eq1.8}
C^{-1}\sup\limits_{Q_{\rho, \theta}(x_0, t_0)} u\leqslant \bigg(\Big[\dfrac{1}{\theta}+\dfrac{1}{\rho^{p^+_{\rho, o}}}\Big]^{1+\frac{N}{p^-_{\rho, o}}}\iint\limits_{Q_{\rho, \theta}(x_0, t_0)} u^{q+1}\,dx dt\bigg)^{\frac{p^-_{\rho, o}}{p^-_{\rho, o}(q+1+N)-(q+1)N}}+\rho 
\end{multline} is valid.
\end{theorem}
\noindent We now deal with the singular case. In order to formulate our next result let us define the following number
\begin{equation}\label{eq1.9}
\varkappa_s:= s\, p^-_{\rho, o} +N\,(p^-_{\rho, o}-q-1).
\end{equation}

\begin{theorem}\label{th1.2}
Let $u$ be a non-negative, local, weak sub-solution to
\eqref{eq1.1}-\eqref{eq1.3} in $\Omega_T$.
Fix $(x_0,t_0)\in\Omega_T$ such that
$Q_{\rho,\theta}(x_0,t_0)\subset
Q_{\rho,\rho}(x_0,t_0)\subset
Q_{8\rho,8\rho}(x_0,t_0)\subset\Omega_T$,
and assume that $q+1>p^+_{\rho,o}$.
Suppose that
\begin{equation}\label{eq1.10}
q>p^-_{\rho,o}-1+(q+1)\frac{p^-_{\rho,o}}{N},
\end{equation}
and let $s>q+1$ be any number such that
\begin{equation}\label{eq1.11}
\varkappa_s>0,\quad \text{i.e.}\quad
s>\frac{N}{p^-_{\rho,o}}(q+1-p^-_{\rho,o}).
\end{equation}
If we assume additionally that
\begin{equation}\label{eq1.12}
\iint\limits_{Q_{\rho,\theta}(x_0,t_0)}
u^s\,dx\,dt<\infty,
\end{equation}
then $u$ is bounded in $Q^-_{\rho,\theta}(x_0,t_0)$.
Moreover, the following quantitative estimate holds:
\begin{equation}\label{eq1.13}
C_1^{-1} \sup\limits_{Q_{\rho, \theta}(x_0, t_0)} u \leqslant \,\bigg(\frac{1}{\theta^{1+\frac{N}{p^-_{\rho, o}}}}
\iint\limits_{Q_{\rho, \theta}(x_0, t_0)} u^{s}\,dx dt\bigg)^{\frac{p^-_{\rho, o}}{\varkappa_s}}+ \Big(\dfrac{\theta}{\rho^{p^+_{\rho, o}}}\Big)^{\frac{1}{q+1-p^+_{\rho, o}}} + \rho,
\end{equation}
where $C_1>0$ is a constant depending only on the data.
The same conclusion holds if
\begin{equation}\label{eq1.14}
q=p^-_{\rho,o}-1+(q+1)\frac{p^-_{\rho,o}}{N},
\end{equation}
provided that \eqref{eq1.11} and \eqref{eq1.12}
are satisfied.
\end{theorem}

\section{Preliminaries}
In this section we give the definition of local weak solution to \eqref{eq1.1} in both variational and mollified weak formulations, so that their equivalence can be easily inferred. Then we prove, thanks to the mollified weak formulation, the energy estimates which are required for the De Giorgi iteration lying at the core of Theorems \ref{th1.1}-\ref{th1.2}.

\subsection{Notion of Solution}

A function $u$ is a  non-negative, local, weak  sub (super)-solution to Eq. \eqref{eq1.1} if 
\begin{equation}\label{eq2.1}
u \in  C_{loc}(0, T; L^{1+q}_{loc}(\Omega))\cap  L^{p_2}_{loc}(0, T; W^{1, p_2}_{loc}(\Omega)),
\end{equation}
and for any compact set $E \subset \Omega$
and every subinterval $[t_{1}, t_{2}]\subset (0, T]$ there holds
\begin{equation}\label{eq2.2}
\int\limits_{E}u^q \varphi\, dx \bigg|^{t_{2}}_{t_{1}} + \int\limits^{t_{2}}\limits_{t_{1}}\int\limits_{E}
\{-u^q \varphi_{\tau}+ \mathbf{A}(x, \tau, u, D u) D \varphi\}\, dx d\tau \leqslant (\geqslant)\, 0,
\end{equation}
 for any testing functions $\varphi \in W^{1,1+q}_{loc}(0, T; L^{1+q}(E))\cap L^{p_2}_{loc}(0, T; W^{1, p_2}_{0}(E))$, $\varphi \geqslant 0$.

\begin{remark}
In the case if $p(x, t)$ is independent of $t$, i.e. $p(x, t)\equiv p(x)$ we give more precise definition of the solution. 
We say that a function $u$ is a  non-negative, local, weak  (sub) super-solution to Eq. \eqref{eq1.1} if 
\begin{equation}\label{eq2.3}
u \in  C_{loc}(0, T; L^{1+q}_{loc}(\Omega))\cap  L_{loc}(0, T; W_{loc}(\Omega)),
\end{equation}
and the integral inequalities \eqref{eq2.2} hold
 for any testing functions $\varphi \in W^{1,1+q}_{loc}(0, T; L^{1+q}(E))\cap L_{loc}(0, T; W_{0}(E)), \zeta\geqslant 0$.
 \end{remark}
\noindent A function $u$ that is both a non-negative, local weak sub and  super-solution to  \eqref{eq1.1} is a non-negative local
weak solution.
\begin{remark}
The structure conditions \eqref{eq1.2} and the function spaces \eqref{eq2.1} and \eqref{eq2.3}
 guarantee the convergence of integrals in \eqref{eq2.2}. For what concerns the existence of the solutions we refer the reader to \cite{AntShm}.
\end{remark}

\subsection{Mollification in Time}

The time derivative of a weak solution exists in the sense of distribution only.   In 
order to overcome the lack of regularity in the time variable, we define the following 
mollification in time:
$$\llbracket v \rrbracket_h(x, t):=\frac{1}{h}\int\limits_0^t e^{\frac{\tau-t}{h}}\,v(x, \tau)\,d\tau,\quad v\in L^1(\Omega_T).$$
We refer the reader to \cite{KinLin} for various properties of this mollification.
\begin{lemma}\label{lem2.1}
For any $s\geqslant 1$ we have
\vskip0.2cm \noindent 
$(i)$ If $v\in L^s(\Omega_T)$, then $\llbracket v \rrbracket_h\in L^s(\Omega_T)$. Moreover, $||\llbracket v \rrbracket_h||_{L^s(\Omega_T)}
\leqslant ||v||_{L^s(\Omega_T)}$, and $\llbracket v \rrbracket_h \rightarrow v$ in $L^s(\Omega_T)$ and almost everywhere on $\Omega_T$
as $h\rightarrow 0$.
\vskip0.2cm \noindent 
$(ii)$ $\llbracket v \rrbracket_h \in C(0, T; L^s(\Omega))$.
\vskip0.2cm \noindent 
$(iii)$ Almost everywhere on $\Omega_T$ there hold
$$\frac{\partial}{\partial t} \llbracket v \rrbracket_h =\frac{1}{h}\big(v-\llbracket v \rrbracket_h\big).$$
\vskip0.2cm \noindent 
$(iV)$ If $Dv \in L^s(\Omega_T)$, then $D \llbracket v \rrbracket_h = \llbracket D v \rrbracket_h \rightarrow D v$ in $L^s(\Omega_T)$
and almost everywhere on $\Omega_T$ as $h \rightarrow 0$.
\vskip0.2cm \noindent 
$(V)$ If $v\in C(0, T; L^s(\Omega))$, then $\llbracket v \rrbracket_h(\cdot, t) \rightarrow v(\cdot, t)$ in $L^s(\Omega)$ and almost everywhere on $\Omega$ for every $t\in (0, T)$ as $h\rightarrow 0$.
\end{lemma}

\noindent From the  differential 
inequalities \eqref{eq2.2} we deduce the mollified version (cf. \cite{KinLin})
\begin{equation}\label{eq2.4}
\iint\limits_{\Omega_T}\Big(\frac{\partial}{\partial t} \llbracket u^q \rrbracket_h\,\varphi+\llbracket \mathbf{A}(x, \tau, u, D u) \rrbracket_h D \varphi\Big)\, dx d\tau \leqslant (\geqslant)\, 0,
\end{equation}
for any non-negative $\varphi\in L^{p_2}(0, T; W^{1, p_2}_0(\Omega))$ \big( $\varphi \in L(0, T; W_0(\Omega))$ if $p(x, t)\equiv p(x)$\big).

\begin{remark}
In the case $q=1$, spaces \eqref{eq2.3} can be also used for $p(x, t)$ satisfying general conditions \eqref{eq1.3}. In this case the Steklov averages in $t$ and Friedrichs' mollifiers in $x$ can be used, see \cite{AlkZik} for the details.
\end{remark}
\noindent We briefly show how to infer the convergence of the mollified energy term. Under conditions \eqref{eq2.1}, by the Young and H\"{o}lder inequalities we have
$$|\llbracket D u \rrbracket_h|^{p(x, t)}\leqslant |\llbracket D u \rrbracket_h|^{p_2}+1\leqslant \frac{1}{h}\int\limits^t_0 e^{\frac{\tau-t}{h}} |D u(x, \tau)|^{p_2} d\tau +1.$$
Integrating this inequality with respect to $t$ over $[0,T]$ and reversing the order of
integration in the double integral and then integrating this inequality with respect to $x$, we arrive at
$$\iint\limits_{\Omega_T} |D\llbracket u\rrbracket_h|^{p(x, t)} dx dt\leqslant ||D u||_{L^{p_2}(\Omega_T)}+ |\Omega_T|.$$
From this, by Lebesque theorem and Lemma \ref{lem2.1},  we obtain
\begin{equation}\label{eq2.5}
\lim\limits_{h\rightarrow 0}\iint\limits_{\Omega_T}|D\llbracket u\rrbracket_h - D u|^{p(x, t)}\,dx dt=0.
\end{equation}
Similarly, under conditions \eqref{eq2.3}
$$|\llbracket D u\rrbracket_h|^{p(x)}\leqslant \frac{1}{h}\int\limits^t_0 e^{\frac{\tau-t}{h}} |D u(x, \tau)|^{p(x)} d\tau,$$
and hence,
$$\iint\limits_{\Omega_T}|D\llbracket u\rrbracket_h|^{p(x)}\,dx dt\leqslant \iint\limits_{\Omega_T}| D u|^{p(x)}\,dx dt,$$
which by Lebesque theorem and Lemma \ref{lem2.1} yields
\begin{equation}\label{eq2.6}
\lim\limits_{h\rightarrow 0}\iint\limits_{\Omega_T}|D\llbracket u\rrbracket_h - D u|^{p(x)}\,dx dt=0.
\end{equation}

\subsection{ Local Energy Estimates}

\begin{lemma}\label{lem2.2}
Let $u$ be a non-negative, local, weak sub (super)-solution to \eqref{eq1.1}-\eqref{eq1.3} in $\Omega_T$. Then there exists $\gamma >0$ depending only on the data, such that for every cylinder $Q_{r,\eta}(y, \tau)\subset \Omega_T$, any $k>0$ and any smooth $\zeta(x, t)$ which is zero for $(x, t)\in \partial B_r(y)\times (\tau-\eta, \tau)$ one has
\begin{multline}\label{eq2.7}
\sup\limits_{\tau-\eta \leqslant t\leqslant t}\int\limits_{B_r(y)\times\{t\}} g_{\pm}(u, k)\,\zeta^{p_2}\,dx+
\iint\limits_{Q_{r, \eta}(y, \tau)}|D(u-k)_{\pm}|^{p(x, t)}\,\zeta^{p_2}\,dx dt\leqslant
\int\limits_{B_r(y)\times\{\tau-\eta\}} g_{\pm}(u, k)\,\zeta^{p_2}\,dx\\+\gamma \iint\limits_{Q_{r, \eta}(y, \tau)}
\Big[g_{\pm}(u, k)|\zeta_t|\zeta^{p_2-1}+(u-k)_{\pm}^{p(x, t)}|D \zeta|^{p(x, t)}\zeta^{p_2-p(x,t)}\big]\,dx dt.
\end{multline}
Here
$$g_{\pm}(u, k):=q\,\int\limits_0^{(u-k)_{\pm}}(z \pm k)^{q-1}\,z\,dz.$$
\end{lemma}
\begin{proof}
We only consider the case of a local weak sub-solution, the case of super-solution can be considered completely similarly. For
fixed $\tau-\eta<t_1<t_2<\tau$ and sufficiently small $\epsilon \in (0,1)$ define
\begin{equation*}
\psi_\epsilon(t):=
\begin{cases}
0,\quad \tau-\eta\leqslant t\leqslant t_1-\epsilon,\\
1+\frac{t-t_1}{\epsilon},\quad t_1-\epsilon\leqslant t \leqslant t_1,\\
1,\quad t_1\leqslant t \leqslant t_2,\\
1-\frac{t-t_2}{\epsilon},\quad t_2\leqslant t\leqslant t_2+\epsilon,\\
0,\quad t_2+\epsilon \leqslant t\leqslant \tau,
\end{cases}
\end{equation*}
and choose in \eqref{eq2.4} $\varphi= (u-k)_+\,\zeta^{p_2}(x, t)\,\psi_\epsilon(t).$ For 
the integral in \eqref{eq2.4} containing the time derivative we compute
\begin{align}\label{eq2.8}
  \iint\limits_{\Omega_T}\,\frac{\partial}{\partial t} \llbracket u^q \rrbracket_h\,\varphi\,dx dt&=\iint\limits_{\Omega_T}
\frac{\partial}{\partial t} \llbracket u^q \rrbracket_h\,(u-k)_+\,\zeta^{p_2}\,\psi_\epsilon \,dx dt=
\iint\limits_{\Omega_T}\,\frac{\partial}{\partial t} \llbracket u^q \rrbracket_h\,\big(\llbracket u^q\rrbracket^{\frac{1}{q}}_h -k\big)_+\,
\zeta^{p_2}\,\psi_\epsilon dx dt\notag\\&+\iint\limits_{\Omega_T}\,\frac{\partial}{\partial t} \llbracket u^q \rrbracket_h\,\Big[(u-k)_+-\big(\llbracket u^q\rrbracket^{\frac{1}{q}}_h -k\big)_+\Big]\,
\zeta^{p_2}\,\psi_\epsilon dx dt= \mathcal{I}_1+\mathcal{I}_2.  
\end{align}
To estimate $\mathcal{I}_1$ we note that
\begin{equation*}
\begin{aligned}
\mathcal{I}_1&=\iint\limits_{\Omega_T}\,\frac{\partial}{\partial t}\bigg[\int\limits^{\llbracket u^q\rrbracket_h}_{k^q}(z^{\frac{1}{q}}-k)_+\,dz \,\,\zeta^{p_2}\,\psi_\epsilon\bigg] dxdt-p_2\iint\limits_{\Omega_T}\int\limits^{\llbracket u^q\rrbracket_h}_{k^q}(z^{\frac{1}{q}}-k)_+\,dz 
\,\,\zeta^{p_2-1}\zeta_t\,\psi_\epsilon dx dt\\&-\iint\limits_{\Omega_T}\int\limits^{\llbracket u^q\rrbracket_h}_{k^q}(z^{\frac{1}{q}}-k)_+\,dz 
\,\,\zeta^{p_2}\,\psi'_\epsilon dx dt\geqslant -p_2\int\limits^{t_2+\epsilon}_{t_1-\epsilon}\int\limits_{B_r(y)}\int\limits^{\llbracket u^q\rrbracket_h}_{k^q}(z^{\frac{1}{q}}-k)_+\,dz 
\,\,\zeta^{p_2-1}|\zeta_t|\,dx dt\\&-\int\limits^{t_2+\epsilon}_{t_1-\epsilon}\int\limits_{B_r(y)}\int\limits^{\llbracket u^q\rrbracket_h}_{k^q}(z^{\frac{1}{q}}-k)_+\,dz 
\,\,\zeta^{p_2}\,\psi'_\epsilon dx dt.
\end{aligned}
\end{equation*}
Using $(iii)$ from Lemma \ref{lem2.1}
\begin{equation*}
\mathcal{I}_2=\frac{1}{h}\iint\limits_{\Omega_T}\,\Big[u^q-\llbracket u^q\rrbracket_h\Big]\Big[(u-k)_+-\big(\llbracket u^q\rrbracket^{\frac{1}{q}}_h -k\big)_+\Big]\,
\zeta^{p_2}\,\psi_\epsilon dx dt\geqslant 0.
\end{equation*}
We can pass $h\rightarrow 0$ on the right-hand side of \eqref{eq2.8}, using the fact that $g_{+}(u, k)=\int\limits^{u^q}_{k^q}(z^{\frac{1}{q}}-k)_+ dz$, we obtain
\begin{equation*}
\lim\limits_{h\rightarrow 0}\iint\limits_{\Omega_T}\,\frac{\partial}{\partial t} \llbracket u^q \rrbracket_h\,\varphi\,dx dt\geqslant
-p_2\int\limits^{t_2+\epsilon}_{t_1-\epsilon}\int\limits_{B_r(y)}g_+(u, k) 
\,\,\zeta^{p_2-1}|\zeta_t|\, dx dt-\int\limits^{t_2+\epsilon}_{t_1-\epsilon}\int\limits_{B_r(y)} g_+(u, k)
\,\,\zeta^{p_2}\,\psi'_\epsilon dx dt.
\end{equation*}
Now we can pass to the limit $\epsilon \rightarrow 0$ in this inequality and obtain
\begin{equation*}
\lim\limits_{\epsilon \rightarrow 0}\lim\limits_{h\rightarrow 0}\iint\limits_{\Omega_T}\,\frac{\partial}{\partial t} \llbracket u^q \rrbracket_h\,\varphi\,dx dt\geqslant\int\limits_{B_r(y)} g_+(u, k)
\,\,\zeta^{p_2}\, dx\bigg|^{t_2}_{t_1}-p_2\int\limits^{t_2}_{t_1}\int\limits_{B_r(y)}g_+(u, k) 
\,\,\zeta^{p_2-1}|\zeta_t|\, dx dt.
\end{equation*}
To estimate the diffusion term we pass to the limit $h\rightarrow 0$, use \eqref{eq2.5}, \eqref{eq1.2} and the Young inequality. 
Passing to the limit $\epsilon \rightarrow 0$ and combining the proceeding inequalities, we arrive at
\begin{multline*}
\int\limits_{B_r(y)\times\{t_2\}} g_{+}(u, k)\,\zeta^{p_2}\,dx+
\frac{1}{\gamma}\int\limits^{t_2}_{t_1}\int\limits_{B_{r}(y)}|D(u-k)_{+}|^{p(x, t)}\,\zeta^{p_2}\,dx dt\leqslant
\int\limits_{B_r(y)\times\{t_1\}} g_{+}(u, k)\,\zeta^{p_2}\,dx\\+\gamma \int\limits^{t_2}_{t_1}\int\limits_{B_{r}(y)}
\Big[g_{+}(u, k)|\zeta_t|\zeta^{p_2-1}+(u-k)_{+}^{p(x, t)}|D \zeta|^{p(x, t)}\zeta^{p_2-p(x,t)}\big]\,dx dt.
\end{multline*}
To complete the proof of \eqref{eq2.7}, we first pass in the last inequality to the limit $t_1\rightarrow \tau-\eta$.
Here, we discard the second integral on the left-hand side and take then the essential 
supremum with respect to $t_2\in (\tau-\eta, \tau)$. This leads to an estimate of the essential 
supremum of the first integral. On the other hand, discarding the first integral and 
passing to the limit $t_2\rightarrow \tau$ we deduce a similar estimate for the second integral of the 
left-hand side. Together, this leads to \eqref{eq2.7}, which completes the proof of the lemma.
\end{proof}
\subsection{Auxiliary Lemmas}
In this subsection, we collect several auxiliary results that will be used in the sequel.\\
The following lemma provides a useful algebraic estimate and can be found in \cite[Lemma~3.2]{Bog}.
\begin{lemma}\label{lem2.3}
Let $q > 0$ and $g_{\pm}$ be defined in Lemma \ref{lem2.2}. There exists a constant $\gamma=\gamma(q)>0$
such that, for all $a, b \in \mathbb{R}$, the following inequality holds true:
$$\frac{1}{\gamma}\big(|a|+|b|\big)^{q-1}(a-b)^2_{\pm}\leqslant g_{\pm}(a, b)\leqslant \gamma \big(|a|+|b|\big)^{q-1}(a-b)^2_{\pm}.$$
\end{lemma}
 The next lemma on parabolic Sobolev embedding can be found in \cite[Proposition~3.1, p.~7]{Bened}.
\begin{lemma}\label{lem2.4} 
Let $E$ be a bounded domain in $\mathbb{R}^N$ and $E_T=E\times (0, T]$. Suppose that 
$u \in L^{\infty}(0, T; L^{m}_{\mathrm{loc}}(\Omega))\cap L^{p}(0, T; W^{1, p}_0(\Omega))$ 
for $m, p\geqslant 1$. Setting $p^*=p\frac{N+m}{N}$, we have 
$$ 
\iint\limits_{E_T} u^{p^*} \, dx \, dt
\leqslant
\gamma \Big(\esssup\limits_{t\in (0, T)}\int\limits_{E}u^m\,dx\Big)^{\frac{ p}{ N}}
\iint\limits_{E_T}|D u|^{p}\,dx \, dt
$$
with a positive constant $\gamma$ depending on $N, p$ and $p^*$, but independent of $u$.
\end{lemma}
The following lemma on geometric convergence is a standard tool in iteration procedures (see \cite[Lemma 4.1, p.12]{Bened}).  
\begin{lemma}\label{lem2.5}
Let $\{y_j\}_{j\in \mathbb{N}}$ be a sequence of positive numbers satisfying the recursive inequalities
$$y_{j+1}\leqslant A\,B^j\,y^{1+\delta}_j,\quad j=0, 1, \dots,$$
where $A, B>1$ and $\delta>0$ are given constants. If 
$$y_0\leqslant A^{-\frac{1}{\delta}}\,B^{-\frac{1}{\delta^2}},$$
then $\lim\limits_{j\rightarrow \infty} y_j=0.$
\end{lemma}

\section{Local Boundedness, Proof of Theorems \ref{th1.1} and \ref{th1.2}}
\noindent Let us fix a point $(x_0, t_0)\in \Omega_T$ and parameters $\theta, \rho>0$ such that 
\[
Q_{\rho, \theta}(x_0, t_0)\subset Q_{\rho, \rho}(x_0, t_0)\subset \Omega_T.
\]
Let $k>0$ be a number to be determined. For $j=0,1,\dots$ we define the levels
\[
k_j=k-\frac{k}{2^{j+1}}, \qquad 
\bar{k}_j:=\frac{1}{2}(k_j+k_{j+1}), \qquad 
k'_j:=\frac{1}{2}(\bar{k}_j+k_{j+1}),
\]
and note that $k_j < \bar{k}_j < k'_j < k_{j+1}.$
We also introduce the decreasing sequences of radii and parameters
\[
r_j:= \frac{1}{2}\rho \Big(1+\frac{1}{2^{j}}\Big), \quad 
\eta_j:=\frac{1}{2}\theta \Big(1 +\frac{1}{2^{j}}\Big), \quad 
B_j:=B_{r_j}(x_0),
\]
and the corresponding cylinders
\[Q_j:=B_j\times I_j, \qquad I_j:= (t_0-\eta_j, t_0).\] 
Define the cutoff function $\zeta_j:=\zeta^{(1)}_j(x) \zeta^{(2)}_j(t)$, where $\zeta^{(1)}_j(x)\in C^1_0(B_j)$ and $\zeta^{(2)}_j(t)\in C^1(\mathbb{R})$ satisfy the following conditions:
$$0\leqslant \zeta^{(1)}_j(x)\leqslant 1, \quad \zeta^{(1)}_j(x)=1 \quad \text{in}\quad \quad B_{j+1},\quad \text{and}\quad |D \zeta^{(1)}_{j}(x)|\leqslant \gamma \frac{2^j}{\rho},$$
$$
 0\leqslant \zeta^{(2)}_j(t)\leqslant 1,\quad \zeta^{(2)}_j(t)=
\begin{cases}
0, & t\leqslant t_0-\eta_j,\\
1, & t\geqslant t_0-\eta_{j+1},
\end{cases} \qquad
|\frac{d}{dt}\zeta^{(2)}_j(t)|\leqslant \gamma \frac{2^j}{\theta}.
$$
Now, we observe that on the set $\{u\geqslant k_{j+1}\}$,
\[
1\leqslant \frac{u+k'_j}{u-k'_j}
\leqslant \frac{2u}{u-k'_j}
\leqslant \frac{2k_{j+1}}{k_{j+1}-k'_j}
\leqslant \gamma 2^{j\gamma},
\]
since the function $s\mapsto \frac{2s}{s-k'_j}$
is decreasing for $s>k'_j$.
Consequently, by Lemma \ref{lem2.3}, we have 
\begin{equation}\label{est_g1}
 g_+(u,k'_j)\geqslant \frac{1}{\gamma}
\left(\frac{u+k'_j}{u-k'_j}\right)^{q-1}
(u-k'_j)^{q+1}
\geqslant \frac{1}{\gamma 2^{j\gamma}}
(u-k_{j+1})^{q+1}, \quad \text{if}\quad u\geqslant k_{j+1}.
\end{equation}
Moreover, arguing similarly on the set $\{u\geqslant k'_j\}$,
we get
\[
1\leqslant \frac{u+k'_j}{u-\bar{k}_j}
\leqslant \frac{2u}{u-\bar{k}_j}
\leqslant \frac{2k'_j}{k'_j-\bar{k}_j}
\leqslant \gamma 2^{j\gamma}.
\]
Thus, applying Lemma \ref{lem2.3}, we obtain
\begin{equation}\label{est_g2}
g_+(u, k^{'}_j)\leqslant \gamma 2^{j\gamma}(u-\bar{k}_j)^{1+q}_+.
\end{equation}
Since $p^-_{\rho,o}\leqslant p(x,t)$ and
$a^{p^-_{\rho,o}}\leqslant 1+a^{p(x,t)}$
for all $a\geqslant 0$, we have
\begin{equation*}
\iint\limits_{Q_j}|D (u-k_{j+1})_+|^{p^-_{\rho, o}}\,\zeta^q_j\,dx dt
\leqslant \gamma \iint\limits_{Q_j}|D (u-k_{j+1})_+|^{p(x, t)}\,\zeta^q_j\,dx dt+
\gamma |Q_j\cap\big\{u\geqslant k_{j+1}\big\}|.
\end{equation*}
Setting
\[
\mathcal{P}=(1+q)\vee p^+_{\rho,o},
\]
we will also use the following simple estimate:
\begin{equation}\label{cheb_ineq}
|Q_j\cap\{u\geqslant k_{j+1}\}|
\leqslant \gamma 2^{j\mathcal{P}}k^{-\mathcal{P}}
\iint\limits_{Q_j}(u-k_j)_+^{\mathcal{P}}\,dx\,dt,
\end{equation}
which is obtained by integrating the inequality
\[
1=\left(\frac{k_{j+1}-k_j}{k_{j+1}-k_j}\right)^{\mathcal{P}}
\leqslant
\left(\frac{u-k_j}{k_{j+1}-k_j}\right)^{\mathcal{P}}
\]
over the set $Q_j\cap\{u\geqslant k_{j+1}\}$.

Using Lemma \ref{lem2.2} and estimates \eqref{est_g1}, \eqref{est_g2}, we obtain
\begin{multline}\label{eq3.1}
\sup\limits_{t\in I_j}\int\limits_{B_j}(u-k_{j+1})^{1+q}_+\,\zeta^{p_2}_j\,dx +\iint\limits_{Q_j}|D (u-k_{j+1})_+|^{p^-_{\rho, o}}\,\zeta^{p_2}_j\,dx dt\\\leqslant \frac{\gamma 2^{j\gamma}}{\theta}\iint\limits_{Q_{j}}(u-\bar{k}_{j})^{1+q}_{+}\,\,dx dt+
\gamma 2^{j\gamma}\,\iint\limits_{Q_{j}}\Big(\frac{(u-\bar{k}_{j})_{+}}{\rho}\Big)^{p(x, t)}\,dx\,dt+ \gamma |Q_j\cap\big\{u\geqslant k_{j+1}\big\}|.
\end{multline}
To estimate the second term on the right-hand side of the previous inequality, we note that the inequality
$\Big(\frac{(u-\bar{k}_{j})_{+}}{\rho}\Big)^{p(x, t)} \leqslant \Big(\frac{(u-\bar{k}_{j})_{+}}{\rho}\Big)^{p^+_{\rho, o}}$
fails in the range $0< \frac{(u-\bar{k}_{j})_{+}}{\rho}<1$. 
Therefore, using \eqref{cheb_ineq},
we proceed as follows:
\begin{multline*}
\sup\limits_{t\in I_j}\int\limits_{B_j}(u-k_{j+1})^{1+q}_+\,\zeta^{p_2}_j\,dx +\iint\limits_{Q_j}|D (u-k_{j+1})_+|^{p^-_{\rho, o}}\,\zeta^{p_2}_j\,dx dt
\\\leqslant\frac{\gamma 2^{j\gamma}}{\theta}\iint\limits_{Q_{j}}(u-\bar{k}_{j})^{1+q}_{+}\,\,dx dt+
\gamma 2^{j\gamma}\,\iint\limits_{Q_{j}}\Big[\Big(\frac{(u-\bar{k}_{j})_{+}}{\rho}\Big)^{p^+_{\rho, o}}+\chi_{\{u>\bar{k}_j\}}\Big]\,dx\,dt+ \gamma |Q_j\cap\big\{u\geqslant k_{j+1}\big\}|\\
\leqslant
\gamma 2^{j\gamma}\bigg[\frac{k^{1+q-\mathcal{P}}}{\theta}+\frac{k^{p^+_{\rho, o}-\mathcal{P}}}{\rho^{p^+_{\rho, o}}}+  \frac{1}{k^{\mathcal{P}}}\bigg]\,\iint\limits_{Q_{j}}(u- k_{j})_{+}^{\mathcal{P}}\,dx\,dt.
\end{multline*}
Here we used the inequality
$(u-\bar{k}_j)_+ \geqslant k_{j+1}-\bar{k}_j
= \gamma k2^{-j}$ on $\{u \geqslant k_{j+1}\}$,
together with $(u-\bar{k}_j)_+ \leqslant (u-k_j)_+$.

 The terms in brackets are estimated by stipulating to take
\begin{equation}\label{eq3.2}
k\geqslant
\begin{cases}
\rho+\Big(\dfrac{\theta}{\rho^{p^+_{\rho, o}}}\Big)^{\frac{1}{q+1-p^+_{\rho,o}}},
& \text{if } p^+_{\rho,o}<q+1,\\
\rho+\Big(\dfrac{\rho^{p^+_{\rho, o}}}{\theta}\Big)^{\frac{1}{p^+_{\rho,o}-q-1}},
& \text{if } p^+_{\rho,o}>q+1,\\
\rho,
& \text{if } p^+_{\rho,o}=q+1.
\end{cases}
\end{equation}
Accordingly, we set
\begin{equation}\label{eq3.3}
H(\theta,\rho)=
\begin{cases}
\dfrac{1}{\theta},
& \text{if } p^+_{\rho,o}<q+1,\\
\dfrac{1}{\rho^{p^+_{\rho, o}}},
& \text{if } p^+_{\rho,o}>q+1,\\
\dfrac{1}{\theta}+\dfrac{1}{\rho^{p^+_{\rho, o}}},
& \text{if } p^+_{\rho,o}=q+1.
\end{cases}
\end{equation}
Then, combining \eqref{eq3.1}, \eqref{eq3.2}, \eqref{eq3.3},
we obtain
\begin{equation}\label{eq3.4}
\begin{aligned}
    \sup\limits_{t\in I_j}
\int\limits_{B_j}(u-k_{j+1})^{1+q}_+
\,\zeta^{p_2}_j\,dx
+\iint\limits_{Q_j}|D(u-k_{j+1})_+|^{p^-_{\rho,o}}
\,\zeta^{p_2}_j\,dx\,dt\\
\leqslant
\gamma 2^{j\gamma}\,H(\theta,\rho)
\iint\limits_{Q_j}(u-k_j)_+^{\mathcal{P}}
\,dx\,dt.
\end{aligned}
\end{equation}

\subsection{Proof of Theorem \ref{th1.1}}

By the log-continuity condition \eqref{eq1.3}
on the exponent $p(x,t)$, we can always assume that
\begin{equation}\label{cond_p^+_p^-}
p^+_{\rho,o}<p^-_{\rho,o}
+(q+1)\frac{p^-_{\rho,o}}{N},
\end{equation}
provided that $\rho>0$ is sufficiently small.
Indeed, since $p^-_{\rho,o}$ depends on $\rho$,
we use the fixed lower bound $p_1$, which is
independent of $\rho$, to choose $\rho$ such that
\[
\frac{K_3}{\log\frac{1}{2\rho}}
\leqslant \frac{1}{2}(q+1)\frac{p_1}{N},
\quad \text{i.e.,}\quad
\rho\leqslant \frac{1}{2}
e^{-\frac{2K_3N}{(q+1)p_1}}.
\]
Then, applying the log-continuity condition
\eqref{eq1.3} and recalling that
$p_1\leqslant p^-_{\rho,o}$, we obtain
\[
p^+_{\rho,o}-p^-_{\rho,o}
\leqslant \frac{K_3}{\log\frac{1}{2\rho}}
\leqslant \frac{1}{2}(q+1)\frac{p_1}{N}
< (q+1)\frac{p^-_{\rho,o}}{N}.
\]

\subsubsection{Case $q<p^-_{\rho, o}-1+(q+1)\frac{p^-_{\rho, o}}{N}$}\label{1case}

Since we are considering the case
$q < p^-_{\rho,o}-1+(q+1)\frac{p^-_{\rho,o}}{N}$,
combining this inequality with \eqref{cond_p^+_p^-},
which holds for sufficiently small $\rho$, we obtain
\[
\mathcal{P}<p^-_{\rho,o}
+(q+1)\frac{p^-_{\rho,o}}{N}=:p_*.
\]
Therefore, applying H\"{o}lder's inequality and
Lemma~\ref{lem2.4}, we obtain
 \begin{multline}\label{eq3.5a}
 \begin{aligned}
  y_{j+1} &:= \iint\limits_{Q_{j+1}}(u-k_{j+1})^{\mathcal{P}}_+\,dx dt
\leqslant \iint\limits_{Q_j}\big[(u-k_{j+1})_+\zeta^{p_2}_j\big]^{\mathcal{P}} dx dt\\
&\leqslant \Big(\iint\limits_{Q_j} \big[(u-k_{j+1})_+\zeta^{p_2}_j\big]^{p_*} dx dt\Big)^{\frac{\mathcal{P}}{p_*}}
|Q_j\cap \big\{u\geqslant k_{j+1}\big\}|^{1-\frac{\mathcal{P}}{p_*}}\\
&\leqslant \gamma \Big(\sup_{t\in I_j}\int_{B_j}(u-k_{j+1})^{1+q}_+ \zeta^{p_2}_j\,dx\Big)^{\frac{\mathcal{P} p^-_{\rho, o}}{p_* N}}
\Big(\iint\limits_{Q_j}|D\big((u-k_{j+1})_+\zeta^{p_2}_j\big)|^{p^-_{\rho,o}} dx dt\Big)^{\frac{\mathcal{P} }{p_*}}\\
&\times|Q_j\cap \big\{u\geqslant k_{j+1}\big\}|^{1-\frac{\mathcal{P}}{p_*}}.  
\end{aligned}
\end{multline}
Using \eqref{eq3.4}, we  estimate the gradient term
\begin{multline}\label{eq3.5e}
 \begin{aligned}
&\qquad\iint\limits_{Q_j}|D\big((u-k_{j+1})_+\zeta^{p_2}_j\big)|^{p^-_{\rho,o}} dx dt\\
&\qquad\le \gamma\iint\limits_{Q_j}\Big(
|D(u-k_{j+1})_+|^{p^-_{\rho,o}} \zeta^{p_2 p^-_{\rho,o}}_j 
+ (u-k_{j+1})_+^{p^-_{\rho,o}} |D\zeta_j|^{p^-_{\rho,o}} \zeta^{p^-_{\rho,o}(p_2-1)}_j\Big) dx dt\\
&\qquad\le \gamma 2^{j\gamma} H(\theta,\rho)\iint\limits_{Q_j}(u-k_j)_+^{\mathcal{P}} dx dt
+ \gamma 2^{j\gamma}\frac{1}{\rho^{p^-_{\rho,o}}}
\iint\limits_{Q_j}(u-k_{j+1})_+^{p^-_{\rho,o}} dx dt.
\end{aligned}
\end{multline}
Since $p^-_{\rho,o}\leqslant\mathcal{P}$, using the
elementary comparison of powers and the ordering
of the levels, we obtain
\[
\gamma 2^{j\gamma}\frac{1}{\rho^{p^-_{\rho,o}}}\iint\limits_{Q_j}(u-k_{j+1})_+^{p^-_{\rho,o}}  dx dt
\le \gamma 2^{j\gamma}\frac{1}{\rho^{p^-_{\rho,o}}}\iint\limits_{Q_j}(u-k_j)_+^{\mathcal{P}} dx dt + \gamma 2^{j\gamma}\frac{1}{\rho^{p^-_{\rho,o}}}|Q_j\cap \big\{u\geqslant k_{j+1}\big\}|.
\]
Substituting the previous inequality into \eqref{eq3.5e} and using \eqref{cheb_ineq} on $Q_j\cap \{u > k_{j+1}\}$, we get
\[
\iint\limits_{Q_j}|D\big((u-k_{j+1})_+\zeta^{p_2}_j\big)|^{p^-_{\rho,o}} dx dt
\le  \gamma 2^{j\gamma} H(\theta,\rho)\iint\limits_{Q_j}(u-k_j)_+^{\mathcal{P}} dx dt.
\]
Combining previous inequality, \eqref{eq3.5a}, \eqref{eq3.4} and \eqref{cheb_ineq}, we derive
\begin{equation}\label{eq3.5_new}
y_{j+1}
\le \gamma 2^{j\gamma} k^{-\mathcal{P}(1-\frac{\mathcal{P}}{p_*})}\big[H(\theta, \rho)\big]^{\frac{\mathcal{P}}{p_*}(1+\frac{p^-_{\rho, o}}{N})}
\,y_j^{1+\frac{\mathcal{P} p^-_{\rho, o}}{p_* N}},
\end{equation}
Applying Lemma \ref{lem2.5}, we deduce that $\lim\limits_{j\rightarrow \infty} y_j=0$, provided that $k$ is chosen to satisfy
\begin{equation*}\label{eq3.6}
k^{q+1+N-\mathcal{P}\frac{N}{p^-_{\rho, o}}}\geqslant \gamma \big[H(\theta, \rho)]^{1+\frac{N}{p^-_{\rho, o}}} \iint\limits_{Q_{\rho, \eta}}(u-\frac{k}{2})^{\mathcal{P}}_+ dx dt.
\end{equation*}
From this, taking into account \eqref{eq3.2} we arrive at the required \eqref{eq1.6}-\eqref{eq1.8}.

\subsubsection{Case  $q=p^-_{\rho, o}-1+(q+1)\frac{p^-_{\rho, o}}{N}$}
Set 
$$y_j:=\iint\limits_{Q_j}(u - k_{j})^{p_*}_{+}\,\,dx dt.$$
Using estimate \eqref{eq3.4} and arguing as in the previous  subsection \ref{1case}, we obtain 
\begin{equation}\label{eq3.7}
\begin{aligned}
    y_{j+1}&\leqslant \iint\limits_{Q_j}\big[(u-k_{j+1})_+\zeta^{p_2}_j\big]^{p_*} dx dt\leqslant
\gamma \Big(\sup\limits_{t\in I_j}\int\limits_{B_j}(u-k_{j+1})^{1+q}_+\,\zeta^{p_2}_j\,dx\Big)^{\frac{ p^-_{\rho, o}}{ N}}\\&\times
\Big(\iint\limits_{Q_j}|D [(u-k_{j+1})_+\zeta^{p_2}_j]|^{p^-_{\rho, o}}\,dx dt\Big)\leqslant\gamma 2^{j\gamma}\,\big[H(\theta, \rho)\big]^{1+\frac{p^-_{\rho, o}}{N}}\, y^{1+\frac{p^-_{\rho, o}}{N}}_j.
\end{aligned}
\end{equation}
Hence, by a Lemma \ref{lem2.5}, it follows that
\[\lim\limits_{j\rightarrow \infty} y_j=\iint\limits_{Q_{\frac{\rho}{2},\frac{\theta}{2}}}(u-k)^{p_*}_{+}\,\,dx dt=0,\]
provided that the initial value 
 is sufficiently small, namely
\begin{equation}\label{eq3.8}
y_0=\iint\limits_{Q_{\rho,\theta}}(u- \frac{k}{2})^{p_*}_{+}\,\,dx dt\leqslant \frac{1}{\gamma}\, \big[H(\theta, \rho)\big]^{-1-\frac{N}{p^-_{\rho, o}}}.
\end{equation}
We now verify that this smallness condition can be achieved by a suitable choice of 
$k$. Indeed, observing that
$$\iint\limits_{Q_{\rho,\theta}}(u- \frac{k}{2})^{p_*}_{+}\,\,dx dt\leqslant \iint\limits_{Q_{\rho,\theta}\cap\{u\geqslant \frac{1}{2} k\}}u^{p_*}\,dx dt.$$
Since $u\in L^{p^*}(Q_{\rho,\theta}),$ by the absolute continuity of the Lebesgue integral we obtain
$$\lim\limits_{k\rightarrow \infty}\iint\limits_{Q_{\rho,\theta}\cap\{u\geqslant \frac{1}{2} k\}}u^{p_*}\,dx dt=0,$$
we choose $k$ large enough such that \eqref{eq3.8} holds, and hence 
$$\sup\limits_{Q_{\frac{1}{2}\rho, \frac{1}{2}\theta}(x_0, t_0)} u\leqslant k,$$
this completes the proof of Theorem \ref{th1.1}.

\subsection{Proof of Theorem \ref{th1.2}}

\subsubsection{Case $q>p^-_{\rho, o}-1+(q+1)\frac{p^-_{\rho, o}}{N}$}

Since $q+1>p_*=p^-_{\rho,o}+(q+1)\frac{p^-_{\rho,o}}{N}$,
we cannot proceed as in the proof of Theorem \ref{th1.1}.
Instead, we estimate the integral on the right-hand side of
\eqref{eq3.4}, taking into account that
$\mathcal{P}=(q+1)\vee p^+_{\rho,o}=q+1$. For $s>q+1$, set
\[
\epsilon:=\frac{s-q-1}{s-p_*}\in(0,1)
\]
and
\[
d:=\bigg(\iint\limits_{Q_{\rho,\theta}(x_0,t_0)}
u^s\,dx\,dt\bigg)^{\frac{1}{s}}.
\]
By H\"older's inequality, taking into account that
$\dfrac{q+1-\epsilon p_*}{1-\epsilon}=s$,
we obtain
\begin{multline*}
\iint\limits_{Q_j}(u-k_{j})^{q+1}_+\,dx dt\leqslant \iint\limits_{Q_j}(u-k_{j})^{\epsilon p_*}_+\,u^{q+1-\epsilon p_*}\,dx dt\leqslant
d^{s(1-\epsilon)}\,\bigg(\iint\limits_{Q_j}(u-k_{j})^{ p_*}_+\,dx dt\bigg)^\epsilon.
\end{multline*}
Combining this with \eqref{eq3.4}, we get
\begin{multline}\label{eq3.9}
\sup\limits_{t\in I_j}\int\limits_{B_j}(u-k_{j+1})^{1+q}_+\,\zeta^{p_2}_j\,dx +\iint\limits_{Q_j}|D (u-k_{j+1})_+|^{p^-_{\rho, o}}\,\zeta^{p_2}_j\,dx dt\\\leqslant \gamma 2^{j \gamma}\,H(\theta, \rho)\,d^{s(1-\epsilon)}\,\bigg(\iint\limits_{Q_{j}}(u- k_{j})_{+}^{p_*}\,dx\,dt\bigg)^{\epsilon}
\end{multline}
Similarly to the derivation of \eqref{eq3.7}, applying
Lemma \ref{lem2.4} and \eqref{eq3.9}, we derive
\begin{equation}\label{eq3.10}
y_{j+1}:=\iint\limits_{Q_{j+1}}(u- k_{j+1})_{+}^{p_*}\,dx\,dt\leqslant \gamma 2^{j \gamma}\,\big[H(\theta, \rho)\,d^{s(1-\epsilon)}\big]^{1+\frac{p^-_{\rho, o}}{N}}\,y^{\epsilon\,\frac{N+p^-_{\rho, o}}{N}}_j.
\end{equation}
A direct computation yields
\begin{equation}\label{eq3.11}
\epsilon(N+p^-_{\rho,o})-N
=\frac{(s-q-1)(N+p^-_{\rho,o})-N(s-p_*)}{s-p_*}
=\frac{\varkappa_s}{s-p_*}>0.
\end{equation}
Consequently,
$\epsilon\frac{N+p^-_{\rho,o}}{N}>1$.
Thus, applying Lemma \ref{lem2.5} to
\eqref{eq3.10}, we obtain
$\lim\limits_{j\to\infty}y_j=0$, provided that
\begin{equation}\label{eq3.12}
y_0\leqslant
\frac{1}{\gamma}
\big[H(\theta,\rho)\,d^{s(1-\epsilon)}\big]^{
-\frac{N+p^-_{\rho,o}}
{\epsilon(N+p^-_{\rho,o})-N}}.
\end{equation}
On the other hand, by H\"older's inequality, we have 
\begin{equation*}
y_0=\iint\limits_{Q_{\rho, \theta}(x_0, t_0)} (u-\frac{k}{2})^{p_*}_+\,dx dt\leqslant \iint\limits_{Q_{\rho, \theta}(x_0, t_0)\cap\{u\geqslant \frac{1}{2} k\}}\,u^{p_*}\,dx dt\leqslant d^{p_*}|Q_{\rho, \theta}(x_0, t_0)\cap\{u\geqslant \frac{1}{2} k\}|^{1-\frac{p_*}{s}}.
\end{equation*}
By the same argument as in \eqref{cheb_ineq}, we obtain
\begin{equation*}
|Q_{\rho,\theta}(x_0,t_0)\cap
\{u\geqslant \tfrac{1}{2}k\}|
\leqslant \Big(\frac{2d}{k}\Big)^s.
\end{equation*}
Consequently,
\begin{equation}\label{eq3.13}
y_0\leqslant
\frac{2^{s-p_*}\,d^s}{k^{s-p_*}}.
\end{equation}
Therefore, the condition \eqref{eq3.12} is satisfied
if $k$ is chosen such that
\begin{equation}\label{eq3.14}
\begin{aligned}
k^{s-p_*}
&\geqslant
\gamma\,d^{
s+\frac{s(1-\epsilon)(N+p^-_{\rho,o})}
{\epsilon(N+p^-_{\rho,o})-N}}
[H(\theta,\rho)]^{
\frac{N+p^-_{\rho,o}}
{\epsilon(N+p^-_{\rho,o})-N}}\\
&=
\gamma\,d^{
\frac{sp^-_{\rho,o}}
{\epsilon(N+p^-_{\rho,o})-N}}
[H(\theta,\rho)]^{
\frac{N+p^-_{\rho,o}}
{\epsilon(N+p^-_{\rho,o})-N}}.
\end{aligned}
\end{equation}
Taking into account \eqref{eq3.11}, condition \eqref{eq3.14}
can be rewritten as
\[
k\geqslant
\gamma [H(\theta,\rho)]^{
\frac{N+p^-_{\rho,o}}{\varkappa_s}}
\bigg(
\iint\limits_{Q_{\rho,\theta}(x_0,t_0)}
u^s\,dx\,dt
\bigg)^{
\frac{p^-_{\rho,o}}{\varkappa_s}}.
\]
Consequently, choosing $k$ sufficiently large
to satisfy both \eqref{eq3.2} and \eqref{eq3.14},
we obtain
\[
\sup\limits_{Q_{\frac{1}{2}\rho,
\frac{1}{2}\theta}(x_0,t_0)}u
\leqslant k.
\]
Since $q+1>p^+_{\rho,o}$, we have
$H(\theta,\rho)=\theta^{-1}$, and hence the preceding estimate implies \eqref{eq1.13}.

\subsubsection{Case $q=p^-_{\rho, o}-1+(q+1)\frac{p^-_{\rho, o}}{N}$}

Set
\[
y_j:=\iint\limits_{Q_j}(u-k_j)^{p_*}_+\,dx\,dt.
\]
By \eqref{eq3.7} and Lemma \ref{lem2.5}, we have
$\lim\limits_{j\to\infty}y_j=0$, provided that
\begin{equation}\label{eq3.15}
y_0\leqslant
\frac{1}{\gamma}
[H(\theta,\rho)]^{-1-\frac{N}{p^-_{\rho,o}}},
\end{equation}
where $H(\theta,\rho)$ is defined in \eqref{eq3.3}.
On the other hand, by \eqref{eq3.13}, we have
\[
y_0\leqslant
\frac{2^{s-p_*}d^s}{k^{s-p_*}}.
\]
Therefore, condition \eqref{eq3.15} is satisfied
if $k$ is chosen such that
\[
k^{s-p_*}\geqslant
\gamma[H(\theta,\rho)]^{1+\frac{N}{p^-_{\rho,o}}}
\iint\limits_{Q_{\rho,\theta}(x_0,t_0)}
u^s\,dx\,dt.
\]
Since $p_*=q+1$ in the present case,
a direct computation gives
\[
s-p_*=\frac{\varkappa_s}{p^-_{\rho,o}}.
\]
Consequently, taking into account \eqref{eq3.2}
and \eqref{eq3.3}, we obtain the required estimate
\eqref{eq1.13}. This completes the proof of
Theorem \ref{th1.2}.

\section*{Acknowledgements} Simone Ciani acknowledges the partial funding of GNAMPA (INdAM), the Italian funding PNR (MUR) and the department of Mathematics of the University of Bolonga. Igor
Skrypnik and Mariia Savchenko were supported by the National Research Foundation
of Ukraine, project No. 2025.07/0369 "Qualitative methods of nonlinear analysis of heterogeneous structures".

\bigskip

CONTACT INFORMATION

\medskip
\textbf{Simone Ciani} \\ University of Bologna Alma Mater, Department of Mathematics, Piazza Porta San Donato 5, 40126, Bologna. E-mail: simone.ciani3@unibo.it

\medskip
\textbf{Igor I.~Skrypnik}\\Institute of Applied Mathematics and Mechanics of the
National Academy of Sciences of Ukraine, Batiuka Str. 19, 84116, Sloviansk, Ukraine. E-mail: ihor.skrypnik@gmail.com

\medskip
\textbf{Mariia Savchenko}\\Institute of Applied Mathematics and Mechanics of the
National Academy of Sciences of Ukraine, Batiuka Str. 19, 84116, Sloviansk, Ukraine. E-mail: shan\_maria@ukr.net

\medskip
\textbf{Yevhen Zozulia}\\National University of  Water and Environmental Engineering,  
 Soborna Str. 11,  33028, Rivne, Ukraine. E-mail: albelgen27@gmail.com

\end{document}